# A goodness-of-fit test for parametric and semi-parametric models in multiresponse regression

SONG XI CHEN[1,2] and INGRID VAN KEILEGOM[3]

[1]*Department of Business Statistics and Econometrics, Guanghua School of Management, Peking University, Beijing 100871, China*
[2]*Department of Statistics, Iowa State University, Ames, Iowa 50011-1210, USA.*
*E-mail: songchen@iastate.edu*
[3]*Institute of Statistics, Université catholique de Louvain, Voie du Roman Pays 20, B-1348 Louvain-la-Neuve, Belgium. E-mail: ingrid.vankeilegom@uclouvain.be*

We propose an empirical likelihood test that is able to test the goodness of fit of a class of parametric and semi-parametric multiresponse regression models. The class includes as special cases fully parametric models; semi-parametric models, like the multiindex and the partially linear models; and models with shape constraints. Another feature of the test is that it allows both the response variable and the covariate be multivariate, which means that multiple regression curves can be tested simultaneously. The test also allows the presence of infinite-dimensional nuisance functions in the model to be tested. It is shown that the empirical likelihood test statistic is asymptotically normally distributed under certain mild conditions and permits a wild bootstrap calibration. Despite the large size of the class of models to be considered, the empirical likelihood test enjoys good power properties against departures from a hypothesized model within the class.

*Keywords:* additive regression; bootstrap; empirical likelihood; goodness of fit; infinite-dimensional parameter; kernel estimation; monotone regression; partially linear regression

## 1. Introduction

Suppose $\{(X_i, Y_i)\}_{i=1}^n$ is an independent and identically distributed random vector, where $Y_i$ is a $k$-variate response and $X_i$ a $d$-variate covariate. Let $m(x) = E(Y_i|X_i = x) = (m_1(x), \ldots, m_k(x))$ be the conditional mean consisting of $k$ regression curves on $R^d$ and $\Sigma(x) = \text{Var}(Y_i|X_i = x)$ be a $k \times k$ matrix whose values change along with the covariate.

Let $m(\cdot) = m(\cdot, \theta, g) = (m_1(\cdot, \theta, g), \ldots, m_k(\cdot, \theta, g))$ be a working regression model of which one would like to check its validity. The form of $m$ is known up to a finite-dimensional parameter $\theta$ and an infinite-dimensional nuisance parameter $g$. The model $m(\cdot, \theta, g)$ includes a wide range of parametric and semi-parametric regression models as special cases. In the absence of $g$, the model degenerates to a fully parametric model $m(\cdot) = m(\cdot, \theta)$, whereas the presence of $g$ covers a range of semi-parametric models including the single or multiindex models and partially linear single-index models. The class also includes models with qualitative constraints, like additive models and models with shape constraints. The variable selection problem, the comparison of regression curves and models for the variance function can be covered by the class of $m(\cdot, \theta, g)$ as well.

Multiresponse regression is frequently encountered in applications. In compartment analysis arising in biological and medical studies as well as chemical kinetics (Atkinson and Bogacka (2002)), a multivariate variable is described by a system of differential equations whose solutions satisfy multiresponse regression (Jacquez (1996)). Surface designs and multivariate random vectors are collected as responses of some controlled variables (covariates) of certain statistical experiments. Khuri (2001) proposed using the generalized linear models for modeling such data and Uciński and Bogacka (2005) studied the issue of optimal designs with an objective for discrimination between two multiresponse system models. The monographs by Bates and Watts ((1988), Chapter 4) and Seber and Wild ((1989), Chapter 11) contain more examples of multiresponse regression as well as their parametric inference.

The need for testing multiple curves occurs even in the context of univariate responses $Y_i$. Consider the following heteroscedastic regression model:

$$Y_i = r(X_i) + \sigma(X_i)e_i,$$

where the $e_i$'s are unit residuals such that $E(e_i|X_i) = 0$ and $E(e_i^2|X_i) = 1$, and $r(\cdot)$ and $\sigma^2(\cdot)$ are, respectively, the conditional mean and variance functions. Suppose $r(x, \theta, g)$ and $\sigma^2(x, \theta, g)$ are certain working parametric or semi-parametric models. In this case, the bivariate response vector is $(Y_i, Y_i^2)^{\mathrm{T}}$ and the bivariate model specification $m(x, \theta, g) = (r(x, \theta, g), \sigma^2(x, \theta, g) + r^2(x, \theta, g))^{\mathrm{T}}$.

The aim of the paper is to develop a nonparametric goodness-of-fit test for the hypothesis

$$H_0 : m(\cdot) = m(\cdot, \theta, g), \tag{1.1}$$

for some known $k$-variate function $m(\cdot, \theta, g)$, some finite-dimensional parameter $\theta \in \Theta \subset R^p$ ($p \geq 1$) and some function $g \in \mathcal{G}$ that is a complete metric space consisting of functions from $R^d$ to $R^q$ ($q \geq 1$). We will use two pieces of nonparametric statistical hardware, the kernel regression estimation technique and the empirical likelihood technique, to formulate a test for $H_0$.

In the case of a single regression curve (i.e., $k = 1$), the nonparametric kernel approach has been widely used to construct goodness-of-fit tests for the conditional mean or variance function. Eubank and Spiegelman (1990), Eubank and Hart (1992), Härdle and Mammen (1993), Hjellvik and Tjøstheim (1995), Fan and Li (1996), Hart (1997) and Hjellvik, Yao and Tjøstheim (1998) have developed consistent tests for a parametric

*Empirical likelihood goodness-of-fit test in multiresponse regression* 957specification by employing the kernel smoothing method based on a fixed bandwidth. Horowitz and Spokoiny (2001) propose a test based on a set of smoothing bandwidths in the construction of the kernel estimator. Its extensions are considered in Chen and Gao (2007) for time series regression models and in Rodríguez-Póo, Sperlich and Vieu (2009) for semi-parametric regression models. Other related references can be found in the books by Hart (1997) and Fan and Yao (2003).

The empirical likelihood (EL) (Owen (1988, 1990)) is a technique that allows the construction of a nonparametric likelihood for a parameter of interest in a nonparametric or semi-parametric setting. Despite that it is intrinsically nonparametric, it possesses two important properties of a parametric likelihood: the Wilks' theorem and the Bartlett correction. Qin and Lawless (1994) establish EL for parameters defined by estimating equations, which is the widest framework for EL formulation. Zhang and Gijbels (2003) considered an EL procedure based on a sieve approach, whereas Chen and Cui (2006) show that the EL admits a Bartlett correction under this general framework. Hjort, McKeague and Van Keilegom (2009) consider the properties of the EL in the presence of both finite and infinite-dimensional nuisance parameters as well as when the data dimension is high. See Owen (2001) for a comprehensive overview of the EL method and references therein.

Goodness-of-fit tests based on the EL have been proposed in the literature, which include Li (2003) and Li and Van Keilegom (2002) for survival data; Einmahl and McKeague (2003) for testing some characteristics of a distribution function; and Chen, Härdle and Li (2003) for conditional mean functions with dependent data. Fan and Zhang (2004) propose a sieve EL test for testing a varying-coefficient regression model that extends the generalized likelihood ratio test of Fan, Zhang and Zhang (2001). They demonstrate that the 'Wilks phenomenon' continues to hold under general error distributions. Tripathi and Kitamura (2003) propose an EL test for conditional moment restrictions. The above two tests and the test we are to propose display an interesting diversity in test statistic formulations via the EL. The basic idea of the EL is to maximize an objective function that is a product of probability weights allocated to observations under certain constraints that characterize the functional object to be tested. Fan and Zhang (2004) apply kernel smoothing in both the objective function and the constraints, whereas Tripathi and Kitamura (2003) smooth only the objective function and we will smooth only the constraints. Fan and Zhang's and our test statistics are based on first constructing local statistics over a range of fixed points and then summing them up to form the final test statistic. The formulation in Tripathi and Kitamura (2003) consists of one step with a global objective function over the range of the entire sample. A common feature among the three formulations is that the test statistics are all asymptotically pivotal (Wilks phenomenon). This is due to the EL's ability to internally studentize a statistic via its optimization procedure.

We consider in the present paper tests for a set of multiple semi-parametric regression functions simultaneously. Multiple regression curves exist when the response $Y_i$ is genuinely multivariate, or when $Y_i$ is in fact univariate but we are interested in testing the validity of a set of feature curves; for example, the conditional mean and conditional variance at the same time. Empirical likelihood is a natural device to formulate goodness-of-fit statistics to test multiple regression curves. This is due to EL's built-in feature to standardize a goodness-of-fit distance measure between a fully nonparametric estimate of the target functional curves and its hypothesized counterparts. This feature is well connected to the Wilks phenomenon in the sieve empirical likelihood of Fan and Zhang (2004) and the generalized likelihood ratio for additive models in Fan and Jiang (2005). The standardization carried out by the EL implicitly uses the true covariance matrix function, say $V(x)$ of the kernel estimator $m(\cdot)$, to studentize the distance between $\hat{m}(\cdot)$ and the hypothesized model $m(\cdot, \theta, g)$, so that the goodness-of-fit statistic is an integrated Mahalanobis distance between the two sets of multivariate curves $\hat{m}(\cdot)$ and $m(\cdot, \theta, g)$. This is attractive as we avoid estimating $V(x)$, which can be a daunting task when $k$ is larger than 1. When testing multiple regression curves, there is an intrinsic issue regarding how much each component-wise goodness-of-fit measure contributes to the final test statistic. The EL distributes the weights naturally according to $V^{-1}(x)$. And most attractively, this is done without requiring extra steps of estimation since it comes as a by-product of the internal algorithm. This attraction of the empirical likelihood has been discovered in Tripathi and Kitamura (2003), who propose an empirical likelihood test that can be used to test parametric multiple response regression models.

The main contribution of the proposed test is its ability to test a large class of regression models in the presence of both finite- and infinite-dimensional parameters. The class includes as special cases fully parametric models; semi-parametric models, like the multi-index and the partially linear models; and models with shape constraints, like monotone regression models. It is shown that the EL test statistic is asymptotically normally distributed under certain mild conditions and permits a wild bootstrap calibration. Despite the large size of the class of models to be considered by the proposed test, the test enjoys good power properties against departures from a hypothesized model within the class.

The paper is organized as follows. In the next section we introduce some notation and formulate the EL test statistic. Section 3 is concerned with the main asymptotic results, namely the asymptotic distribution of the test statistic both under the null hypothesis and under a local alternative, and the consistency of the bootstrap approximation. In Section 4 we focus on a number of particular models and apply the general results of Section 3 to these models. Simulation results are reported in Section 5. We conclude the paper by giving in Section 6 the assumptions and the proofs of the main results.

## 2. The test statistic

Let $Y_i = (Y_{i1}, \ldots, Y_{ik})^{\mathrm{T}}$ and $m(x) = (m_1(x), \ldots, m_k(x))^{\mathrm{T}}$, where $m_l(x) = E(Y_{il}|X_i = x)$ is the $l$th regression curve for $l = 1, \ldots, k$. Let $\epsilon_i = Y_i - m(X_i)$ be the $i$th residual vector. Define $\sigma_{lj}(x) = \mathrm{Cov}(\epsilon_{il}, \epsilon_{ij}|X_i = x)$, which is the conditional covariance between the $l$th and $j$th component of the residual vector. Then, the conditional covariance matrix $\Sigma(x) = \mathrm{Var}(Y_i|X_i = x) = (\sigma_{lj}(x))_{k \times k}$.

Let $K$ be a $d$-dimensional kernel with a compact support on $[-1, 1]^d$. Without loss of generality, $K$ is assumed to be a product kernel based on a univariate kernel $k$, that is,

$K(t_1,\ldots,t_d) = \prod_{i=1}^{d} k(t_i)$, where $k$ is a $r$th-order kernel supported on $[-1,1]$ and

$$\int k(u)\,du = 1, \quad \int u^l k(u)\,du = 0 \quad \text{for } l = 1,\ldots,r-1 \quad \text{and} \quad \int u^r k(u)\,du = k_r \neq 0$$

for an integer $r \geq 2$. Define $K_h(u) = h^{-d} K(u/h)$. The Nadaraya–Watson (NW) kernel estimator of $m_l(x)$, $l = 1,\ldots,k$, is

$$\hat{m}_l(x) = \frac{\sum_{i=1}^{n} K_{h_l}(x - X_i) Y_{il}}{\sum_{t=1}^{n} K_{h_l}(x - X_t)},$$

where $h_l$ is the smoothing bandwidth for curve $l$. Different bandwidths are allowed to smooth different curves, which is sensible for multivariate responses. Then

$$\hat{m}(x) = (\hat{m}_1(x),\ldots,\hat{m}_k(x))^{\mathrm{T}}$$

is the kernel estimator of the multiple regression curves. We assume throughout the paper that $h_l/h \to \beta_l$ as $n \to \infty$, where $h$ represents a baseline level of the smoothing bandwidth and $c_0 \leq \min_l\{\beta_l\} \leq \max_l\{\beta_l\} \leq c_1$ for finite and positive constants $c_0$ and $c_1$ free of $n$.

Under the null hypothesis (1.1),

$$Y_i = m(X_i, \theta_0, g_0) + \epsilon_i, \qquad i = 1,\ldots,n, \tag{2.1}$$

where $\theta_0$ is the true value of $\theta$ in $\Theta$, $g_0$ is the true function in $\mathcal{G}$, the errors $\epsilon_1,\ldots,\epsilon_n$ are independent and identically distributed, $E(\epsilon_i | X_i = x) = 0$ and $\mathrm{Var}(\epsilon_i | X_i = x) = \Sigma(x)$.

Let $\hat{\theta}$ be a $\sqrt{n}$-consistent estimator of $\theta_0$ and $\hat{g}$ be a consistent estimator of $g_0$ under a norm $\|\cdot\|_{\mathcal{G}}$ defined on the complete metric space $\mathcal{G}$. Any $\sqrt{n}$-consistent estimator of $\theta_0$ would be fine; for instance, the pseudo-likelihood estimator that assumes that the residual distribution is normal. We suppose $\hat{g}$ is a kernel estimator based on a kernel $L$ of order $s \geq 2$ and a bandwidth sequence $b$, most likely different from the bandwidth $h$ used to estimate $m$. We will require that $\hat{g}$ converges to $g_0$ faster than $(nh^d)^{-1/2}$, the optimal rate in a completely $d$-dimensional nonparametric model. As demonstrated in Section 4, this can be easily satisfied since $g$ is of lower-dimensional than the saturated nonparametric model for $m$.

Each $m_l(x, \hat{\theta}, \hat{g})$ is smoothed by the same kernel $K$ and bandwidth $h_l$ as in the kernel estimator $\hat{m}_l(x)$, in order to prevent the bias of the kernel regression estimators entering the asymptotic distribution of the test statistic (see also Härdle and Mammen (1993)):

$$\tilde{m}_l(x, \hat{\theta}, \hat{g}) = \frac{\sum_{i=1}^{n} K_{h_l}(x - X_i) m_l(X_i, \hat{\theta}, \hat{g})}{\sum_{t=1}^{n} K_{h_l}(x - X_t)}$$

for $l = 1,\ldots,k$. Let $\tilde{m}(x, \hat{\theta}, \hat{g}) = (\tilde{m}_1(x, \hat{\theta}, \hat{g}),\ldots,\tilde{m}_k(x, \hat{\theta}, \hat{g}))^{\mathrm{T}}$.

We note in passing that the dimension of the response $Y_i$ does not contribute to the curse of dimensionality. Rather, it is the dimension of the covariate $X_i$ that contributes, since $X_i$ is the direct target of smoothing. Hence, as far as the curse of dimensionality is concerned, testing multiple curves is the same as testing a single regression curve.

To formulate the empirical likelihood ratio test statistics, we first consider a fixed $x \in R^d$ and then globalize by integrating the local likelihood ratio over a compact set $S \subset R^d$ in the support of $X$. For each fixed $x \in S$, let

$$\hat{Q}_i(x,\hat{\theta}) = (K_{h_1}(x-X_i)(Y_{i1} - \tilde{m}_1(x,\hat{\theta},\hat{g})), \ldots, K_{h_k}(x-X_i)(Y_{ik} - \tilde{m}_k(x,\hat{\theta},\hat{g})))^{\mathrm{T}}, \quad (2.2)$$

which is a vector of local residuals at $x$ and its mean is approximately zero.

Let $\{p_i(x)\}_{i=1}^n$ be non-negative real numbers representing empirical likelihood weights allocated to $\{(X_i, Y_i)\}_{i=1}^n$. The minus 2 log empirical likelihood ratio for the multiple conditional mean evaluated at $\tilde{m}(x,\hat{\theta},\hat{g})$ is

$$\ell\{\tilde{m}(x,\hat{\theta},\hat{g})\} = -2\sum_{i=1}^n \log\{np_i(x)\}$$

subject to $p_i(x) \geq 0$, $\sum_{i=1}^n p_i(x) = 1$ and $\sum_{i=1}^n p_i(x)\hat{Q}_i(x,\hat{\theta}) = 0$. By introducing a vector of Lagrange multipliers $\lambda(x) \in R^k$, a standard empirical likelihood derivation (Owen (1990)) shows that the optimal weights are given by

$$p_i(x) = \frac{1}{n}\{1 + \lambda^{\mathrm{T}}(x)\hat{Q}_i(x,\hat{\theta})\}^{-1}, \quad (2.3)$$

where $\lambda(x)$ solves

$$\sum_{i=1}^n \frac{\hat{Q}_i(x,\hat{\theta})}{1 + \lambda^{\mathrm{T}}(x)\hat{Q}_i(x,\hat{\theta})} = 0. \quad (2.4)$$

Integrating $\ell\{\tilde{m}(x,\hat{\theta},\hat{g})\}$ against a weight function $\pi$ supported on $S$, gives

$$\Lambda_n(\vec{h}) = \int \ell\{\tilde{m}(x,\hat{\theta},\hat{g})\}\pi(x)\,\mathrm{d}x,$$

which is our EL test statistic based on the bandwidth vector $\vec{h} = (h_1, \ldots, h_k)^{\mathrm{T}}$.

Define $\hat{\bar{Q}}(x,\hat{\theta}) = n^{-1}\sum_{i=1}^n \hat{Q}_i(x,\hat{\theta})$, let $f(x)$ be the density of $X$, $R(t) = \int K(u)K(tu)\,\mathrm{d}u$ and $V(x) = f(x)(\beta_j^{-d}R(\beta_l/\beta_j)\sigma_{lj}(x))_{k\times k}$. Note that $R(1) = R(K) =: \int K^2(u)\,\mathrm{d}u$ and that $\beta_j^{-d}R(\beta_l/\beta_j) = \beta_l^{-d}R(\beta_j/\beta_l)$ indicating that $V(x)$ is a symmetric matrix.

Derivations given in Section 6 show that

$$\Lambda_n(\vec{h}) = nh^d \int \hat{\bar{Q}}^{\mathrm{T}}(x,\theta_0)V^{-1}(x)\hat{\bar{Q}}(x,\theta_0)\pi(x)\,\mathrm{d}x + \mathrm{o}_p(h^{d/2}),$$

where $h^{d/2}$ is the stochastic order of the first term on the right-hand side if $d < 4r$. Here $r$ is the order of the kernel $K$. Since $\hat{\bar{Q}}(x,\theta_0) = f(x)\{\hat{m}(x) - \tilde{m}(x,\theta_0,\hat{g})\}\{1 + \mathrm{o}_p(1)\}$, $\hat{\bar{Q}}(x,\theta_0)$ serves as a raw discrepancy measure between $\hat{m}(x) = (\hat{m}_1(x), \ldots, \hat{m}_k(x))$ and the hypothesized model $m(x,\theta_0,\hat{g})$. There is a key issue on how much each $\hat{m}_l(x) - \tilde{m}_l(x,\theta_0,\hat{g})$ contributes to the final statistic. The EL distributes the contributions according to $nh^d V^{-1}(x)$, the inverse of the covariance matrix of $\hat{\vec{Q}}(x, \theta_0)$, which is the most natural choice. The nice thing about the EL formulation is that this is done without explicit estimation of $V(x)$ due to its internal standardization. Estimating $V(x)$ when $k$ is large can be challenging if not just tedious.

## 3. Main results

Let $(\gamma_{lj}(x))_{k \times k} = ((\beta_j^{-d} R(\beta_l/\beta_j) \sigma_{lj}(x))_{k \times k})^{-1}$,

$$\omega_{l_1,l_2,j_1,j_2}(\beta, K) = \int \int \int \beta_{l_2}^{-d} K(u) K(v) K\{(\beta_{j_2} z + \beta_{l_1} u)/\beta_{l_2}\} K(z + \beta_{j_1} v/\beta_{j_2}) \, du \, dv \, dz,$$

$$\sigma^2(K, \Sigma) = 2 \sum_{l_1,l_2,j_1,j_2}^{k} \beta_{l_2}^{-d} \omega_{l_1,l_2,j_1,j_2}(\beta, K) \int \gamma_{l_1 j_1}(x) \gamma_{l_2 j_2}(x) \sigma_{l_1 l_2}(x) \sigma_{j_1 j_2}(x) \pi^2(x) \, dx,$$

which is a bounded quantity under assumptions (A.1) and (A.4) given in Section 6.

**Theorem 3.1.** *Under the assumptions* (A.1)–(A.6) *and* (B.1)–(B.5) *given in Section 6, and under* $H_0$,

$$h^{-d/2}\{\Lambda_n(\vec{h}) - k\} \xrightarrow{d} N(0, \sigma^2(K, \Sigma))$$

*as* $n \to \infty$.

*Remark 3.1 (Equal bandwidths).* If $h_1 = \cdots = h_k = h$, that is, $\beta_1 = \cdots = \beta_k = 1$, then $\omega_{l_1,l_2,j_1,j_2}(\beta, K) = K^{(4)}(0)$, where $K^{(4)}$ is the convolution of $K^{(2)}$ and $K^{(2)}$ is the convolution of $K$, that is, $K^{(2)}(u) = \int K(v) K(u + v) \, dv$. Since $V(x) = f(x) R(K) \Sigma(x)$ in the case of equal bandwidths, $\sum_{l=1}^{k} \gamma_{lj_1} \sigma_{lj_2}(x) = I(j_1 = j_2) R^{-1}(K)$, where $I$ is the indicator function. Therefore, $\sigma^2(K, \Sigma) = 2k K^{(4)}(0) R^{-2}(K) \int \pi^2(x) \, dx$, which is entirely known upon giving the kernel function $K$ and the weight function $\pi$. Hence, the EL test statistic is asymptotically pivotal.

*Remark 3.2 (Unequal bandwidths).* If the bandwidths are not all the same, the asymptotic variance of $\Lambda_n(\vec{h})$ may depend on $\Sigma(x)$, which means that the EL test statistic is no longer asymptotically pivotal. However, the distribution of $\Lambda_n(\vec{h})$ is always free of the design distribution of $X_i$.

Let $q_{n\alpha}$ be the upper $\alpha$-quantile of the distribution of $h^{-d/2}\{\Lambda_n(\vec{h}) - k\}$ for a significance level $\alpha \in (0, 1)$. Theorem 3.1 assures that $q_{n\alpha} \to z_\alpha$, the upper $\alpha$ quantile of $N(0, \sigma^2(K, \Sigma))$. However, the convergence can be slow. There is also an issue of estimating $\sigma^2(K, \Sigma)$ when different bandwidths are used. For these reasons we prefer to use a bootstrap approximation to calibrate the quantile $q_{n\alpha}$.

**Remark 3.3 (Bootstrap).** Let $\hat{\epsilon}_i = Y_i - \hat{m}(X_i)$ be the estimated residual vectors for $i=1,\ldots,n$ and $G$ be a multivariate $k$-dimensional random vector such that $E(G)=0$, $\text{Var}(G)=I_k$ and $G$ has bounded fourth-order moments. To facilitate simple construction of the test statistic, and faster convergence, we propose the following bootstrap estimate of $q_{n\alpha}$.

Step 1: For $i=1,\ldots,n$, generate $\epsilon_i^* = \hat{\epsilon}_i G_i$ where $G_1,\ldots,G_n$ are independent and identical copies of $G$, and let $Y_i^* = m(X_i,\hat{\theta},\hat{g}) + \epsilon_i^*$. Re-estimate $\theta$ and $g$ based on $\{(X_i,Y_i^*)\}_{i=1}^n$ and denote them as $\hat{\theta}^*$ and $\hat{g}^*$.

Step 2: Compute the EL ratio at $\tilde{m}(x,\hat{\theta}^*,\hat{g}^*)$ based on $\{(X_i,Y_i^*)\}_{i=1}^n$, denote it as $\ell^*\{\tilde{m}(x,\hat{\theta}^*,\hat{g}^*)\}$ and then obtain the bootstrap version of the test statistic $\Lambda_n^*(\vec{h}) = \int \ell^*\{\tilde{m}(x,\hat{\theta}^*,\hat{g}^*)\}\pi(x)\,dx$ and let $\xi^* = h^{-d/2}\{\Lambda_n^*(\vec{h}) - k\}$.

Step 3: Repeat steps 1 and 2 $N$ times, and obtain $\xi_1^* \leq \cdots \leq \xi_N^*$ without loss of generality.

The bootstrap estimate of $q_{n\alpha}$ is then $\hat{q}_{n\alpha} =: \xi_{[N\alpha]+1}^*$.

*The proposed EL test* with $\alpha$-level of significance rejects $H_0$ if $h^{-d/2}\{\Lambda_n(\vec{h}) - k\} > \hat{q}_{n\alpha}$.

**Remark 3.4 (Bandwidth selection).** Each bandwidth $h_l$ used in the kernel regression estimator $\hat{m}_l(x)$ can be chosen by a standard bandwidth selection procedure; for instance, the cross-validation (CV) method. The range in terms of order of magnitude for all the $k$ bandwidths $\{h_l\}_{l=1}^k$ covers the order of $n^{-1/(d+2r)}$, which is the optimal order that minimizes the mean integrated squared error in the estimation of $m_l$ and is also the asymptotic order of the bandwidth selected by the CV method. We also note that once $\{h_l\}_{l=1}^k$ are chosen, the same set of bandwidths will be used in formulating the bootstrap version of the test statistic $\Lambda_n^*(\vec{h})$.

**Theorem 3.2.** *Under assumptions* (A.1)–(A.6) *and* (B.1)–(B.5) *given in Section 6, and under* $H_0$,

$$P(h^{-d/2}\{\Lambda_n(\vec{h}) - k\} \geq \hat{q}_{n\alpha}) \to \alpha,$$

*as* $\min(n,N) \to \infty$.

Theorem 3.2 maintains that the proposed test has asymptotically correct size.

We next consider the power of the test under a sequence of local alternatives. First, consider the following local alternative hypothesis:

$$H_{1n}: m(\cdot) = m(\cdot,\theta_0,g_0) + c_n \Gamma_n(\cdot), \quad (3.1)$$

where $c_n = n^{-1/2}h^{-d/4}$ and $\Gamma_n(x) = (\Gamma_{n1}(x),\ldots,\Gamma_{nk}(x))^\mathrm{T}$ for some bounded functions $\Gamma_{nl}(\cdot)$ ($l=1,\ldots,k$).

We need the following theorem on the asymptotic distribution of the EL test statistic under $H_{1n}$ in order to evaluate the property of the EL test.

**Theorem 3.3.** *Under the assumptions* (A.1)–(A.7) *and* (B.1)–(B.5) *given in Section 6, and under* $H_{1n}$,

$$h^{-d/2}\{\Lambda_n(\vec{h}) - k\} \xrightarrow{d} N(\beta(f, K, \Sigma, \Gamma), \sigma^2(K, \Sigma))$$

*as* $n \to \infty$, *where*

$$\beta(f, K, \Sigma, \Gamma) = \sum_{l,j=1}^{k} \int \Gamma_l(x)\Gamma_j(x)\gamma_{lj}(x)f(x)\pi(x)\,dx$$

$$= \int \Gamma^{\mathrm{T}}(x)V^{-1}(x)\Gamma(x)f^2(x)\pi(x)\,dx$$

*and* $\Gamma(x) = \lim_{n \to \infty} \Gamma_n(x)$, *assuming such a limit does exists.*

**Remark 3.5 (Power).** The asymptotic mean of the EL test statistic is given by $\int \Gamma^{\mathrm{T}}(x)V^{-1}(x)\Gamma(x)f^2(x)\pi(x)\,dx$, which is bounded away from zero since $V(x)$ is positive definite with the smallest eigenfunction uniformly bounded away from zero. As a result, the EL test has a non-trivial asymptotic power,

$$\Phi[\{\beta(f, K, \Sigma, \Gamma) - z_\alpha\}/\sigma(K, \Sigma)],$$

where $\Phi$ is the distribution function of the standard normal distribution. We note here that the above power is attained for any $\Gamma(x)$ without requiring specific directions in which $H_{1n}$ deviates from $H_0$. This indicates the proposed test is able to test consistently any departure from $H_0$. In summary we have the following theorem.

**Remark 3.6 (Choice of $\pi$).** The choice of the weight function $\pi$ will affect the performance of the test. This is reflected in the power function given in Remark 3.5 as both the $\beta$ and the $\sigma$ function depend on $\pi$. One possible way to select $\pi$ is to maximize the following expression with respect to $\pi$:

$$\{\beta(f, K, \Sigma, \Gamma) - z_\alpha\}/\sigma(K, \Sigma),$$

which is the argument of the asymptotic power function of the test. Here both $\beta(f, K, \Sigma, \Gamma)$ and $\sigma(K, \Sigma)$ depend on the choice of $\pi$. The power also depends on the local alternative function $\Gamma_n$ as well as on $f$ and the covariance $\Sigma$. While $f$ and $\Sigma$ can be estimated empirically, the optimization of the selection of $\pi$ will have to be done by assuming some special form of $\Gamma_n$.

**Theorem 3.4.** *Under assumptions* (A.1)–(A.6) *and* (B.1)–(B.5) *given in Section 6, and under* $H_{1n}$,

$$P(h^{-d/2}\{\Lambda_n(\vec{h}) - k\} \geq \hat{q}_{n\alpha}) \to \Phi[\{\beta(f, K, \Sigma, \Gamma) - z_\alpha\}/\sigma(K, \Sigma)],$$

*as* $\min(n, N) \to \infty$.

## 4. Examples

In this section we will apply the general results obtained in Section 3 on a number of particular models: partially linear models, single index models, additive models, monotone regression models, the selection of variables and the simultaneous testing of the conditional mean and variance. These six examples form a representative subset of the more complete list of examples listed in the introduction section. For the other examples not treated here, the development is quite similar.

### 4.1. Partially linear models

Consider the model

$$Y_i = m(X_i, \theta_0, g_0) + \epsilon_i$$
$$= \theta_{00} + \theta_{01}X_{i1} + \cdots + \theta_{0,d-1}X_{i,d-1} + g_0(X_{id}) + \epsilon_i, \tag{4.1}$$

where $Y_i$ is a one-dimensional response variable ($k=1$), $d > 1$, $E(\epsilon_i | X_i = x) = 0$ and $\mathrm{Var}(\epsilon_i | X_i = x) = \Sigma(x)$ ($1 \le i \le n$). For identifiability reasons we assume that $E(g_0(X_{id})) = 0$. This testing problem has been studied in Yatchew (1992), Whang and Andrews (1993) and Rodríguez-Póo, Sperlich and Vieu (2009), among others. For any $\theta \in R^d$ and $x \in R$, let

$$\hat{h}(x, \theta) = \sum_{i=1}^{n} W_{in}(x, b)[Y_i - \theta_0 - \theta_1 X_{i1} - \cdots - \theta_{d-1} X_{i,d-1}], \tag{4.2}$$

$$\hat{g}(x, \theta) = \hat{h}(x, \theta) - \frac{1}{n}\sum_{i=1}^{n} \hat{h}(X_{id}, \theta), \tag{4.3}$$

where

$$W_{in}(x, b) = \frac{L((x - X_{id})/b)}{\sum_{j=1}^{n} L((x - X_{jd})/b)},$$

$b$ is a univariate bandwidth sequence and $L$ is a kernel function. Next, define

$$\hat{\theta} = \underset{\theta \in R^d}{\mathrm{argmin}} \sum_{i=1}^{n} [Y_i - \theta_0 - \theta_1 X_{i1} - \cdots - \theta_{d-1} X_{i,d-1} - \hat{g}(X_{id}, \theta)]^2.$$

Then, $\hat{\theta} - \theta_0 = \mathrm{O}_p(n^{-1/2})$, see Härdle, Liang and Gao ((2000), Chapter 2), and

$$|m(X_i, \theta_0, \hat{g}) - m(X_i, \theta_0, g_0)| = |\hat{g}(X_i, \theta_0) - g_0(X_i)|$$
$$= \mathrm{O}_p\{(nb)^{-1/2}\log(n)\} = \mathrm{o}_p\{(nh^d)^{-1/2}\log(n)\},$$

uniformly in $1 \leq i \leq n$, provided $h^d/b \to 0$. This is the case when $h \sim n^{-1/(d+4)}$ and $b \sim n^{-1/5}$. Hence, condition (B.1) is satisfied. Conditions (B.2) and (B.3) obviously hold, since

$$\frac{\partial m(X_i, \theta_0, g)}{\partial \theta} = (1, X_{i1}, \ldots, X_{i,d-1})^{\mathrm{T}} \quad \text{and} \quad \frac{\partial^2 m(X_i, \theta_0, g)}{\partial \theta \, \partial \theta^{\mathrm{T}}} = 0$$

for any $g$. When the order of the kernel $L$ equals 2, $E\{\hat{g}(x, \theta_0)\} = g_0(x) + \mathrm{O}(b^2)$ uniformly in $x$ and $\mathrm{O}(b^2)$ is $\mathrm{o}(h^2)$ provided $b/h \to 0$, which is satisfied for the above choices of $h$ and $b$. Hence, (B.4) is satisfied for $r = 2$. Finally, for condition (B.5), if we choose $\mathcal{G}$ equal to the class of continuously differentiable and bounded functions, then it is easily seen that $P(\hat{g} \in \mathcal{G}) \to 1$ as $n \to \infty$.

## 4.2. Single index models

In single index models it is assumed that

$$Y_i = m(X_i, \theta_0, g_0) + \epsilon_i = g_0(\theta_0^{\mathrm{T}} X_i) + \epsilon_i, \tag{4.4}$$

where $k$ (the dimension of $Y_i$) equals 1, $\theta_0 = (\theta_{01}, \ldots, \theta_{0d})^{\mathrm{T}}$, $X_i = (X_{i1}, \ldots, X_{id})^{\mathrm{T}}$ for some $d > 1$, $E(\epsilon_i | X_i = x) = 0$ and $\mathrm{Var}(\epsilon_i | X_i = x) = \Sigma(x)$ ($1 \leq i \leq n$). In order to identify the model, set $\|\theta_0\| = 1$. See, e.g., Xia, Li, Tong and Zhang (2004), Stute and Zhu (2005) and Rodríguez-Póo, Sperlich and Vieu (2009) for procedures to test this single index model. For any $\theta \in \Theta$ and $u \in R$, let

$$\hat{g}(u, \theta) = \sum_{i=1}^{n} \frac{L_b(u - \theta^{\mathrm{T}} X_i)}{\sum_{j=1}^{n} L_b(u - \theta^{\mathrm{T}} X_j)} Y_i.$$

Then, the estimator of $\theta_0$ is defined by

$$\hat{\theta} = \operatorname*{argmin}_{\theta: \|\theta\|=1} \sum_{i=1}^{n} [Y_i - \hat{g}(\theta^{\mathrm{T}} X_i, \theta)]^2.$$

Härdle, Hall and Ichimura (1993) showed that $\hat{\theta} - \theta_0 = \mathrm{O}_p(n^{-1/2})$. Obviously, from standard kernel regression theory we know that

$$\max_i |m(X_i, \theta_0, \hat{g}) - m(X_i, \theta_0, g_0)| \leq \sup_u |\hat{g}(u, \theta_0) - g_0(u)|$$

$$= \mathrm{O}_p\{(nb)^{-1/2} \log(n)\} = \mathrm{o}_p\{(nh^d)^{-1/2} \log(n)\},$$

$$\max_i \left| \frac{\partial}{\partial \theta} m(X_i, \theta_0, \hat{g}) - \frac{\partial}{\partial \theta} m(X_i, \theta_0, g_0) \right| \leq C \sup_u |\hat{g}'(u, \theta_0) - g_0'(u)|$$

$$= \mathrm{O}_p\{(nb^3)^{-1/2} \log(n)\} = \mathrm{o}_p(1),$$

$$\max_i \left| \frac{\partial^2}{\partial \theta \, \partial \theta^{\mathrm{T}}} m(X_i, \theta_0, \hat{g}) \right| \leq C \sup_u |\hat{g}''(u, \theta_0)|$$

$$= C \sup_u |g_0''(u)| + \mathrm{O}_p\{(nb^5)^{-1/2} \log(n)\} = \mathrm{o}_p(n^{1/2})$$

and $\sup_u |E\{\hat{g}(u,\theta_0)\} - g_0(u)| = \mathrm{O}(b^2) = \mathrm{o}(h^2)$, for some $C > 0$, provided $h^d/b \to 0$ and $nb^3 \log^{-2}(n) \to \infty$, which is the case (for the partially linear model) when, for example, $h \sim n^{-1/(d+4)}$ and $b \sim n^{-1/5}$.

### 4.3. Additive models

We suppose now that the model is given by

$$Y_i = m_{00} + m_{10}(X_{i1}) + \cdots + m_{d0}(X_{id}) + \epsilon_i, \qquad (4.5)$$

where $k = 1$, $d > 1$, $E(\epsilon_i|X_i = x) = 0$, $\mathrm{Var}(\epsilon_i|X_i = x) = \Sigma(x)$ and $E(m_{j0}(X_{ij})) = 0$ ($1 \le i \le n$; $1 \le j \le d$). The estimation of the parameter $m_{00}$ and of the functions $m_{j0}(\cdot)$ ($1 \le j \le d$) has been considered in Linton and Nielsen (1995) (marginal integration); Opsomer and Ruppert (1997) (backfitting); and Mammen, Linton and Nielsen (1999) (smooth backfitting). Using the covering technique to extend pointwise convergence results to uniform results (see, e.g., Bosq (1998)), it can be shown that the estimators $\hat{m}_j(\cdot)$ ($j = 1, \ldots, d$) considered in these papers satisfy the following properties:

$$\sup_x |\hat{m}_j(x) - m_{j0}(x)| = \mathrm{O}_p\{(nb)^{-1/2} \log(n)\},$$

$$\sup_x |E\{\hat{m}_j(x)\} - m_{j0}(x)| = \mathrm{O}(b^2),$$

where $b$ is the bandwidth used for either of these estimators. Hence, assumptions (B.1)–(B.5) hold true provided $h^d/b \to 0$ and $b/h \to 0$, which is the case when, for example, $h$ and $b$ equal the optimal bandwidths for kernel estimation in dimension $d$, respectively 1, namely $h \sim n^{-1/(d+4)}$ and $b \sim n^{-1/5}$ (take $r = s = 2$).

### 4.4. Monotone regression

Consider now the following model

$$Y_i = m_0(X_i) + \epsilon_i, \qquad (4.6)$$

where $X_i$ and $Y_i$ are one-dimensional and we assume that $m_0$ is monotone. An overview of nonparametric methods for estimating a monotone regression function, as well as testing for monotonicity, is given in Gijbels (2005). See also Dette, Neumeyer and Pilz (2006) for a recent contribution in this area. In the latter paper, a preliminary (not necessarily monotone) estimator of the regression function is used to obtain a monotone estimator of the inverse of the regression function, which is then inverted to obtain a monotone estimator of the regression function itself.

Let $\hat{m}(x)$ be an arbitrary estimator of $m_0(x)$ under the assumption of monotonicity that is based on a bandwidth sequence $b$ and a kernel $L$ of order $s$ and that satisfies

$$\sup_x |\hat{m}(x) - m_0(x)| = O_p\{(nb)^{-1/2}\log(n)\},$$

$$\sup_x |E\{\hat{m}(x)\} - m_0(x)| = O(b^s)$$

(as for the additive model, the uniformity in $x$ can be obtained by using classical tools based on the covering technique). For instance, the estimator given in Dette, Neumeyer and Pilz (2006) satisfies this property (see their Theorem 3.2). Let $\mathcal{G}$ be the class of monotone functions defined on the support of $X$. Then the regularity conditions (B.1)–(B.5) on $\hat{m}(x)$ are satisfied provided $h/b \to 0$ and $b^s/h^r \to 0$; for example, when $s=3$, $r=2$, $b = Kn^{-1/5}$ and $h = b\log^{-1}(n)$. Note that conditions (B.2)–(B.3) are automatically satisfied, since there is no parametric component in the model. Contrary to the previous examples, here we cannot take $h$ and $b$ equal to the optimal bandwidths for kernel estimation, as they both involve univariate smoothing. It now follows from Theorem 3.1 that $h^{-1/2}(\Lambda_n(h) - 1) \xrightarrow{d} N(0, \sigma^2)$ for some $\sigma^2 > 0$, which by Remark 3.1 only depends on $K$ and $\pi$.

## 4.5. Selection of variables

In this example we apply the general testing procedure on the problem of selecting explanatory variables in regression. Let $X_i = (X_i^{(1)\mathrm{T}}, X_i^{(2)\mathrm{T}})^\mathrm{T}$ be a vector of $d = d_1 + d_2$ ($d_1, d_2 \geq 1$) explanatory variables. We like to test whether the vector $X_i^{(2)}$ should or should not be included in the model. See Delgado and González Manteiga (2001) for other nonparametric approaches to this problem. Our null model is

$$Y_i = m_0(X_i^{(1)}) + \epsilon_i. \tag{4.7}$$

Hence, under the hypothesized model, the regression function $m(x^{(1)}, x^{(2)})$ is equal to a function $m_0(x^{(1)})$ only. In our testing procedure we estimate $m_0(\cdot)$ by

$$\hat{m}(x^{(1)}) = \sum_{i=1}^n \frac{L_b(x^{(1)} - X_i^{(1)})}{\sum_j L_b(x^{(1)} - X_j^{(1)})} Y_i,$$

where $L$ is a $d_1$-dimensional kernel function of order $s = 2$ and $b$ a bandwidth sequence. It is easily seen that this estimator satisfies the regularity conditions provided $h^d/b^{d_1} \to 0$ and $b/h \to 0$ (take $r = 2$). As before, the optimal bandwidths for estimation, namely $h \sim n^{-1/(d+4)}$ and $b \sim n^{-1/(d_1+4)}$ satisfy these constraints.

## 4.6. Simultaneous testing of the conditional mean and variance

Let $Z_i = r(X_i) + \Sigma^{1/2}(X_i)e_i$ where $Z_i$ is a $k_1$-dimensional response variable of a $d$-dimensional covariate $X_i$, and $r(x) = E(Z_i|X_i = x)$ and $\Sigma(x) = \mathrm{Var}(Z_i|X_i = x)$ are re-spectively the conditional mean and variance functions. This is a standard multivariate nonparametric regression model. Suppose that $r(x,\theta,g)$ and $\Sigma(x,\theta,g)$ are certain working models for the conditional mean and variance, respectively. Hence, the hypothesized regression model is

$$Z_i = r(X_i,\theta,g) + \Sigma^{1/2}(X_i,\theta,g)e_i, \tag{4.8}$$

where the standardized residuals $\{e_i\}_{i=1}^n$ satisfy $E(e_i|X_i) = 0$ and $\text{Var}(e_i|X_i) = I_d$. Here, $I_d$ is the $d$-dimensional identity matrix. Clearly, the parametric (without $g$) or semiparametric (with $g$) model specification of (4.8) consists of two components of specifications: one for the regression part $r(X_i,\theta,g)$ and the other is the conditional variance part $\Sigma(X_i,\theta,g)$. The model (4.8) is valid if and only if both components of the specifications are valid simultaneously. Hence, we need to test the goodness of fit of both $r(x,\theta,g)$ and $\Sigma(x,\theta,g)$ simultaneously.

To use the notation of this paper, we have

$$m(x,\theta,g) = (r(x,\theta,g), \text{vec}\{\Sigma(x,\theta,g) + r(x,\theta,g)r^{\text{T}}(x,\theta,g)\})^{\text{T}}$$

and the multivariate "response" $Y_i = (Z_i, \text{vec}(Z_i Z_i^{\text{T}}))^{\text{T}}$. Here $\text{vec}(A)$ denotes the operator that stacks columns of a matrix $A$ into a vector.

## 5. Simulations

We carry out two simulation studies. For the first one, consider the following model:

$$Y_i = 1 + 0.5 X_{i1} + a g_1(X_{i1}) + g_2(X_{i2}) + \varepsilon_i \tag{5.1}$$

($i = 1,\ldots,n$). Here, the covariates $X_{i1}$ and $X_{i2}$ are independent and follow a uniform distribution on $[0,1]$, and the error $\varepsilon_i$ is independent of $X_i = (X_{i1}, X_{i2})$ and follows a normal distribution with mean zero and variance given by $\text{Var}(\varepsilon_i|X_i) = (1.5 + X_{i1} + X_{i2})^2/100$. Several choices are considered for the constant $a \geq 0$ and the functions $g_1$ and $g_2$. We are interested in testing whether the data follow a partially linear model, in the sense that the regression function is linear in $X_{i1}$, and (possibly) nonlinear in $X_{i2}$.

We will compare our EL-based test with the test considered by Rodríguez-Póo, Sperlich and Vieu (2009) (RSV hereafter), which is based on the $L_\infty$-distance between a completely nonparametric kernel estimator of the regression function and an estimator obtained under the assumption that the model is partially linear.

The simulations are carried out for samples of size 100 and 200. The significance level is $\alpha = 0.05$. A total of 300 samples are selected at random, and for each sample 300 random resamples are drawn. We choose to work with the weight function $\pi(x) = I(0.1 \leq x \leq 0.9)$, in order to avoid boundary effects for small and large values of $x$. A triangular kernel function $K(u) = (1 - |u|)I(|u| \leq 1)$ is used and we determine the bandwidth $b$ by using a cross-validation procedure. For the bandwidth $h$, we follow the procedure used by Rodríguez-Póo, Sperlich and Vieu (2009), that is, we consider the test statistic $\sup_{h_0 \leq h \leq h_1}[h^{-d/2}\{\Lambda_n(\vec{h}) - k\}]$, where $h_0$ and $h_1$ are chosen in such a way that the bandwidth obtained by cross-validation is included in the interval. For $n = 100$ we take $h_0 = 0.22$ and $h_1 = 0.28$ and for $n = 200$ we select $h_0 = 0.18$ and $h_1 = 0.24$. The critical values for this test statistic are obtained from the distribution of the bootstrap statistic, given by $\sup_{h_0 \leq h \leq h_1}[h^{-d/2}\{\Lambda_n^*(\vec{h}) - k\}]$.

The results are shown in Table 1. The table shows that the level is well respected for both sample sizes, and for both choices of the function $g_2$. Under the alternative hypothesis, all the considered models demonstrate that the power increases with increasing sample size and increasing value of $a$. The empirical likelihood test is in general more powerful than the RSV test when $a$ is small ($a = 0.5$ and $1.0$). For the largest $a$ considered, that is, $a = 1.5$, the RSV test is slightly more powerful. However, this happens when both tests enjoy a large amount of power.

We now consider a second simulation study, based on the following model:

$$Y_i = 1 + 0.5X_{i1} + aX_{i1}^2 + \exp(X_{i2}) + 0.15\exp(cX_{i1})e_i \qquad (5.2)$$

($i = 1, \ldots, n$). As before, the covariates $X_{ij}$ ($j = 1, 2$) are independent and follow a uniform distribution on $[0, 1]$ and the error $\varepsilon_i$ is independent of $X_i = (X_{i1}, X_{i2})$ and follows a standard normal distribution. We are interested in simultaneous testing of the regression and variance function for several choices of $a$ and $c$. The null model corresponds to $a = c = 0$, that is, under the null hypothesis we have a homoscedastic partial linear model. The same choices for $n, \alpha, K, b$ and $h$ are taken as in the first simulation study. As before, we carry out 300 simulations and, for each sample, 300 random resamples are generated. The weight function is now given by $\pi(x) = \pi(x_1, x_2) = \prod_{j=1}^{2} I(0.1 \leq x_j \leq 0.9)$.

The results are shown in Table 2. As far as we know, there is no competitor for this test in the literature. As is clear from the table, the rejection probabilities are close to the nominal level under the null hypothesis and increase when $a$, $c$ and $n$ get larger.

**Table 1.** Rejection probabilities under the null hypothesis ($a = 0$) and under the alternative hypothesis ($a > 0$). The test of Rodríguez-Póo, Sperlich and Vieu (2009) is indicated by 'RSV', the new test is indicated by 'EL' for model (5.1)

| | | $g_1(x_1) = x_1^2$ | | | | $g_1(x_1) = 2\log(x_1 + 0.5)$ | | | |
| | | $n = 100$ | | $n = 200$ | | $n = 100$ | | $n = 200$ | |
| $g_2(x_2)$ | $a$ | RSV | EL | RSV | EL | RSV | EL | RSV | EL |
| --- | --- | --- | --- | --- | --- | --- | --- | --- | --- |
| $\exp(x_2)$ | 0 | 0.047 | 0.053 | 0.040 | 0.043 | 0.047 | 0.053 | 0.040 | 0.043 |
| | 0.5 | 0.123 | 0.153 | 0.160 | 0.193 | 0.123 | 0.147 | 0.127 | 0.160 |
| | 1 | 0.377 | 0.420 | 0.653 | 0.683 | 0.387 | 0.400 | 0.657 | 0.660 |
| | 1.5 | 0.787 | 0.743 | 0.973 | 0.980 | 0.747 | 0.723 | 0.973 | 0.960 |
| $\frac{2}{x_2+1}$ | 0 | 0.033 | 0.037 | 0.043 | 0.053 | 0.033 | 0.037 | 0.043 | 0.053 |
| | 0.5 | 0.110 | 0.120 | 0.153 | 0.177 | 0.107 | 0.133 | 0.113 | 0.147 |
| | 1 | 0.373 | 0.397 | 0.667 | 0.657 | 0.407 | 0.440 | 0.660 | 0.713 |
| | 1.5 | 0.753 | 0.733 | 0.977 | 0.963 | 0.797 | 0.763 | 0.990 | 0.983 |

# 6. Assumptions and proofs

*Assumptions.*

(A.1) $K$ is a $d$-dimensional product kernel of the form $K(t_1,\ldots,t_d) = \prod_{j=1}^{d} k(t_j)$, where $k$ is an $r$th-order ($r \geq 2$) univariate kernel (i.e., $k(t) \geq 0$ and $\int k(t)\,\mathrm{d}t = 1$) supported on $[-1,1]$ and $k$ is symmetric, bounded and Lipschitz continuous.

(A.2) The baseline smoothing bandwidth $h$ satisfies $nh^{d+2r} \to K$ for some $K \geq 0$, $nh^{3d/2}\log^{-4}(n) \to \infty$, and $h_l/h \to \beta_l$ as $n \to \infty$, where $c_0 \leq \min_{1 \leq l \leq k}\{\beta_l\} \leq \max_{1 \leq l \leq k}\{\beta_l\} \leq c_1$ for finite and positive constants $c_0$ and $c_1$. Moreover, $d < 4r$ and the weight function $\pi$ is bounded, Lipschitz continuous on its compact support $S$ and satisfies $\int \pi(x)\,\mathrm{d}x = 1$.

(A.3) Let $\epsilon_i = Y_i - m(X_i, \theta_0, g_0) = (\epsilon_{i1},\ldots,\epsilon_{ik})^{\mathrm{T}}$. $E(|\prod_{j=1}^{6} \epsilon_{il_j}||X_i = x)$ is uniformly bounded for all $l_1,\ldots,l_6 \in \{1,\ldots,k\}$ and all $x \in S$.

(A.4) $f(x)$ and all the $\sigma_{lj}^2(x)$'s have continuous derivatives up to the second order in $S$, $\inf_{x \in S} f(x) > 0$ and $\min_l \inf_{x \in S} \sigma_{ll}^2(x) > 0$. Let $\xi_1(x)$ and $\xi_k(x)$ be the smallest and largest eigenvalues of $V(x)$. We assume that $c_2 \leq \inf_{x \in S} \xi_1(x) \leq \sup_{x \in S} \xi_k(x) \leq c_3$ for finite and positive constants $c_2$ and $c_3$.

(A.5) $\Theta$ is a subspace of $R^p$, $P(\hat{\theta} \in \Theta) \to 1$ as $n \to \infty$, and $\hat{\theta}$ satisfies $\hat{\theta} - \theta_0 = \mathrm{O}_p(n^{-1/2})$.

(A.6) $m(x,\theta,g)$ is twice continuously partially differentiable with respect to the components of $\theta$ and $x$ for all $g$, and the derivatives are bounded uniformly in $x \in S, \theta \in \Theta$ and $g \in \mathcal{G}$.

(A.7) The functions $\Gamma_{nl}(x)$ ($l = 1,\ldots,k$) appearing in the local alternative hypothesis converge to $\Gamma_l(x)$ as $n \to \infty$, and $\Gamma_l(x)$ is uniformly bounded with respect to $x$.

Let $\hat{\Delta}_l(x,\theta) = \tilde{m}_l(x,\theta,\hat{g}) - \tilde{m}_l(x,\theta,g_0)$ for $l = 1,\ldots,k$, $\hat{\Delta}(x,\theta) = (\hat{\Delta}_l(x,\theta))_{l=1}^{k}$,

$$\hat{Q}_i^{(2)}(x,\theta) = (K_{h_1}(x - X_i)\hat{\Delta}_1(x,\theta),\ldots,K_{h_k}(x - X_i)\hat{\Delta}_k(x,\theta))^{\mathrm{T}},$$

and let $\|\cdot\|$ be the Euclidean norm.

The following conditions specify stochastic orders for some quantities involving $\hat{Q}_i^{(2)}(x,\theta_0)$. They can be verified for particular choices of the null model. In Section 4, we

**Table 2.** Rejection probabilities under the null hypothesis ($a = c = 0$) and under the alternative hypothesis ($a, c > 0$) for model (5.2)

| $a$ | $c$ | $n = 100$ | $n = 200$ |
|---|---|---|---|
| 0   | 0 | 0.053 | 0.043 |
| 0.5 | 1 | 0.080 | 0.130 |
| 1   | 2 | 0.203 | 0.310 |
| 2   | 4 | 0.550 | 0.757 |
| 3   | 6 | 0.853 | 0.913 |

*show that these conditions hold true for the common semi-parametric regression models, like the single index and partial linear model.*

(B.1) $\max_{i,l} |m_l(X_i, \theta_0, \hat{g}) - m_l(X_i, \theta_0, g_0)| = o_p\{(nh^d)^{-1/2} \log(n)\}$.
(B.2) $\max_{i,l} |\frac{\partial m_l(X_i, \theta_0, \hat{g})}{\partial \theta} - \frac{\partial m_l(X_i, \theta_0, g_0)}{\partial \theta}| = o_p(1)$.
(B.3) $\max_{i,l} \|\frac{\partial^2 m_l(X_i, \theta_0, \hat{g})}{\partial \theta \partial \theta^{\mathrm{T}}}\| = o_p(n^{1/2})$.
(B.4) $\sup_{x \in S} \|E\{m(x, \theta_0, \hat{g})\} - m(x, \theta_0, g_0)\| = o(h^r)$.
(B.5) $P(\hat{g} \in \mathcal{G}) \to 1$ as $n \to \infty$.

Because of space restrictions, the proofs of the next five lemmas are omitted. They can be found in a technical report, available from the authors (Chen and Van Keilegom (2009)).

**Lemma 6.1.** *Assumption* (B.1) *implies* (B.1a)–(B.1c), *given by*

(B.1a) $\sup_{x \in S}[n^{-1} \sum_{i=1}^n \hat{Q}_i^{(2)}(x, \theta_0) \hat{Q}_i^{(2)\mathrm{T}}(x, \theta_0)] = o_p\{n^{-1} h^{-2d} \log^2(n)\}$.
(B.1b) $\sup_{x \in S} \max_{1 \le i \le n} \|\hat{Q}_i^{(2)}(x, \theta_0)\| = o_p\{n^{-1/2} h^{-3d/2} \log(n)\}$.
(B.1c) $\sup_{x \in S}[n^{-1} \sum_{i=1}^n \hat{Q}_i^{(2)}(x, \theta_0)] = o_p\{(nh^d)^{-1/2} \log(n)\}$.

Next, we consider the uniform rate of convergence of the Lagrange multiplier $\lambda(x)$.

**Lemma 6.2.** *Assume* (A.1)–(A.6) *and* (B.1)–(B.5). *Then, under* $H_0$,

$$\sup_{x \in S} \|\lambda(x) h^{-d}\| = O_p\{(nh^d)^{-1/2} \log(n)\}.$$

The following lemma gives a one-step expansion for $\lambda(x)$.

**Lemma 6.3.** *Assume* (A.1)–(A.6) *and* (B.1)–(B.5). *Then, under* $H_0$,

$$\lambda(x) h^{-d} = V^{-1}(x) \hat{\tilde{Q}}(x, \hat{\theta}) + O_p\{(nh^d)^{-1} \log^3(n)\}, \tag{6.1}$$

*uniformly with respect to* $x \in S$.

We next derive an expansion of the EL ratio statistic.

**Lemma 6.4.** *Assume* (A.1)–(A.6) *and* (B.1)–(B.5). *Then, under* $H_0$,

$$\ell\{\tilde{m}(x, \hat{\theta}, \hat{g})\} = nh^d \hat{\tilde{Q}}^{\mathrm{T}}(x, \hat{\theta}) V^{-1}(x) \hat{\tilde{Q}}(x, \hat{\theta}) + \hat{q}_n(x, \hat{\theta}) + o_p(h^{d/2}),$$

*uniformly with respect to* $x \in S$, *where*

$$\hat{q}_n(x, \hat{\theta}) = nh^d \hat{\tilde{Q}}^{\mathrm{T}}(x, \hat{\theta})\{(h^d \hat{S}_n(x, \hat{\theta}))^{-1} - V^{-1}(x)\}\hat{\tilde{Q}}(x, \hat{\theta}) + \tfrac{2}{3} nh^d \hat{D}_n(x).$$

Applying Lemma 6.4, it can be shown that the EL test statistic can be written as

$$\Lambda_n(\vec{h}) = nh^d \int \hat{\tilde{Q}}^{\mathrm{T}}(x,\theta_0) V^{-1}(x) \hat{\tilde{Q}}(x,\theta_0) \pi(x) \, \mathrm{d}x + R_n + \mathrm{o}_p(h^{d/2}), \qquad (6.2)$$

where

$$R_n = \int \hat{q}_n(x,\hat{\theta}) \pi(x) \, \mathrm{d}x + 2nh^d \int \hat{\tilde{Q}}^{\mathrm{T}}(x,\theta_0) V^{-1}(x) \{\hat{\tilde{Q}}(x,\hat{\theta}) - \hat{\tilde{Q}}(x,\theta_0)\} \pi(x) \, \mathrm{d}x$$

$$+ nh^d \int \{\hat{\tilde{Q}}(x,\hat{\theta}) - \hat{\tilde{Q}}(x,\theta_0)\}^{\mathrm{T}} V^{-1}(x) \{\hat{\tilde{Q}}(x,\hat{\theta}) - \hat{\tilde{Q}}(x,\theta_0)\} \pi(x) \, \mathrm{d}x.$$

Let us consider the order of $R_n$.

**Lemma 6.5.** *Assume* (A.1)–(A.6) *and* (B.1)–(B.5). *Then, under* $H_0$, $R_n = \mathrm{o}_p(h^{d/2})$.

**Lemma 6.6.** *Under assumptions* (A.1)–(A.6) *and* (B.1)–(B.5), *and under* $H_0$,

$$\Lambda_n(\vec{h}) = \Lambda_{n1}(\vec{h}) + \mathrm{o}_p(h^{d/2}), \qquad (6.3)$$

*where* $\Lambda_{n1}(\vec{h}) = nh^d \int \hat{\tilde{Q}}^{(1)\mathrm{T}}(x,\theta_0) V^{-1}(x) \hat{\tilde{Q}}^{(1)}(x,\theta_0) \pi(x) \, \mathrm{d}x$.

**Proof.** Lemma 6.5 and (6.2) lead to

$$\Lambda_n(\vec{h}) = \Lambda_{n1}(\vec{h}) + 2nh^d \int \hat{\tilde{Q}}^{(1)\mathrm{T}}(x,\theta_0) V^{-1}(x) \hat{\tilde{Q}}^{(2)}(x,\theta_0) \pi(x) \, \mathrm{d}x$$

$$+ nh^d \int \hat{\tilde{Q}}^{(2)\mathrm{T}}(x,\theta_0) V^{-1}(x) \hat{\tilde{Q}}^{(2)}(x,\theta_0) \pi(x) \, \mathrm{d}x + \mathrm{o}(h^{d/2}).$$

Applying the same analysis to the term $\hat{D}_{n3}(x)$ in the proof of Lemma 6.5, we have

$$nh^d \int \hat{\tilde{Q}}^{(2)\mathrm{T}}(x,\theta_0) V^{-1}(x) \hat{\tilde{Q}}^{(2)}(x,\theta_0) \pi(x) \, \mathrm{d}x = \mathrm{o}_p\{(nh^d)^{-1} \log^2(n)\} = \mathrm{o}_p(h^{d/2}).$$

It remains to check the order of $\Lambda_{n2}(\vec{h}) = nh^d \int \hat{\tilde{Q}}^{(1)\mathrm{T}}(x,\theta_0) V^{-1}(x) \hat{\tilde{Q}}^{(2)}(x,\theta_0) \pi(x) \, \mathrm{d}x$. Applying the same style of derivation as for $\hat{D}_{n1}(x)$, it can be shown that $\Lambda_{n2}(\vec{h}) = \mathrm{o}_p(h^{d/2})$. This finishes the proof. □

**Proof of Theorem 3.1.** First note that

$$\Lambda_{n1}(\vec{h}) = n^{-1} h^d \sum_{i,j}^{n} \sum_{l,t}^{k} \tilde{\epsilon}_{il}(x) \tilde{\epsilon}_{jt}(x) \int K_{h_l}(x - X_i) K_{h_t}(x - X_j) \gamma_{lt}(x) f^{-1}(x) \pi(x) \, \mathrm{d}x,$$

where $\tilde{\epsilon}_{il}(x) = Y_{il} - \tilde{m}_l(x,\theta_0,g_0)$. Let $K^{(2)}(\beta_l, \beta_t, u) = \beta_t^{-d} \int K(z) K(\frac{\beta_l z}{\beta_t} + u) \, \mathrm{d}z$, which is a generalization of the standard convolution of $K$ to accommodate different bandwidths and is symmetric with respect to $\beta_l$ and $\beta_t$. By a change of variable and noticing that $K$ is a compact kernel supported on $[-1, 1]^d$, $\Lambda_{n1}(\vec{h}) = \Lambda_{n11}(\vec{h})\{1 + \mathrm{O}_p(h^2)\}$, where

$$\begin{aligned}
\Lambda_{n11}(\vec{h}) &= n^{-1} \sum_{i,j}^{n} \sum_{l,t}^{k} \epsilon_{il}\epsilon_{jt} K^{(2)}\left(\beta_l, \beta_t, \frac{X_i - X_j}{h_t}\right) \sqrt{\frac{\pi(X_i)\pi(X_j)\gamma_{lt}(X_i)\gamma_{lt}(X_j)}{f(X_i)f(X_j)}} \\
&= n^{-1} \sum_{i \neq j}^{n} \sum_{l,t}^{k} \epsilon_{il}\epsilon_{jt} K^{(2)}\left(\beta_l, \beta_t, \frac{X_i - X_j}{h_t}\right) \sqrt{\frac{\pi(X_i)\pi(X_j)\gamma_{lt}(X_i)\gamma_{lt}(X_j)}{f(X_i)f(X_j)}} \\
&\quad + n^{-1} \sum_{i=1}^{n} \sum_{l,t}^{k} \epsilon_{il}\epsilon_{it} \beta_t^{-d} R(\beta_l/\beta_t) \frac{\pi(X_i)\gamma_{lt}(X_i)}{f(X_i)} \\
&=: \Lambda_{n111}(\vec{h}) + \Lambda_{n112}(\vec{h}).
\end{aligned}$$
(6.4)

It is straightforward to show that $\Lambda_{n112}(\vec{h}) = k + \mathrm{o}_p(h^{d/2})$. Thus, it contributes only to the mean of the test statistic. As $\Lambda_{n111}(\vec{h})$ is a degenerate $U$-statistic with kernel depending on $n$, straightforward but lengthy calculations lead to

$$h^{-d/2}\Lambda_{n111}(\vec{h}) \xrightarrow{d} N(0, \sigma^2(K, \Sigma)).$$

The establishment of the above asymptotic normality can be achieved by either an approach using the martingale central limit theorem (Hall and Heyde (1980)) as demonstrated in Hall (1984) or an approach using the generalized quadratic forms (de Jong (1987)) as demonstrated in Härdle and Mammen (1993). Note that $(nh^d)^{-1}\log^4(n) = \mathrm{o}(h^{d/2})$. Applying Slutsky's theorem leads to the result. $\square$

**Proof of Theorem 3.2.** It can be checked that given the original sample $\chi_n = \{(X_i, Y_i)\}_{i=1}^n$, versions of assumptions (B.1)–(B.5) are true for the bootstrap resample. Hence Lemmas 6.2–6.6 are valid for the resample given $\chi_n$. In particular, let $\hat{\tilde{Q}}^*(x, \hat{\theta})$ be the bootstrap version of $\hat{\tilde{Q}}(x, \theta_0)$, let $\hat{V}(x) = \hat{f}(x)(\beta_j^{-d} R(\beta_l/\beta_j)\hat{\sigma}_{lj}(x))_{k \times k}$, where $\hat{\sigma}_{lj}(x) = \hat{f}^{-1}(x)n^{-1}\sum_{i=1}^{n} K_h(x - X_i)\hat{\epsilon}_{il}\hat{\epsilon}_{ij}$, $\hat{f}(x) = n^{-1}\sum_{i} K_h(x - X_i)$, and let $(\hat{\gamma}_{lj}(x))_{k \times k} = \hat{f}(x)\hat{V}^{-1}(x)$. Then, conditional on $\chi_n$, $\ell^*\{\tilde{m}(x, \hat{\theta}^*, \hat{g}^*)\} = nh^d \hat{\tilde{Q}}^{*\mathrm{T}}(x, \hat{\theta}) \times \hat{V}^{-1}(x)\hat{\tilde{Q}}^*(x, \hat{\theta}) + \mathrm{o}_p(h^{d/2})$ and $\Lambda_n^*(\vec{h}) = \Lambda_{n11}^*(\vec{h}) + \mathrm{o}_p(h^{d/2})$, where

$$\Lambda_{n11}^*(\vec{h}) = n^{-1} \sum_{i,j}^{n} \sum_{l,t}^{k} \epsilon_{il}^* \epsilon_{jt}^* K^{(2)}\left(\beta_l, \beta_t, \frac{X_i - X_j}{h_t}\right) \sqrt{\frac{\pi(X_i)\pi(X_j)\hat{\gamma}_{lt}(X_i)\hat{\gamma}_{lt}(X_j)}{\hat{f}(X_i)\hat{f}(X_j)}},$$

which are respectively the bootstrap versions of (6.3) and (6.4).

Then apply the central limit theorem for degenerate $U$-statistics as in the proof of Theorem 3.1, conditional on $\chi_n$, $h^{-d/2}(\Lambda_{n11}^* - k) \xrightarrow{d} N(0, \sigma^2(K, \hat{\Sigma}))$, where $\sigma^2(K, \hat{\Sigma})$ is $\sigma^2(K, \Sigma)$ with $\Sigma(x)$ replaced by $\hat{\Sigma}(x) = (\hat{\sigma}_{lj}(x))_{k \times k}$. This implies that

$$h^{-d/2}(\Lambda_n^* - k) \xrightarrow{d} N(0, \sigma^2(K, \hat{\Sigma})). \tag{6.5}$$

Let $\hat{Z} \stackrel{d}{=} N(0, \sigma^2(K, \hat{\Sigma}))$ and $Z \stackrel{d}{=} N(0, \sigma^2(K, \Sigma))$, and $\hat{z}_\alpha$ and $z_\alpha$ be the upper-$\alpha$ quantiles of $N(0, \sigma^2(K, \hat{\Sigma}))$ and $N(0, \sigma^2(K, \Sigma))$, respectively. Recall that $\hat{q}_{n\alpha}$ and $q_{n\alpha}$ are, respectively, the upper-$\alpha$ quantile of $h^{-d/2}(\Lambda_n^* - k)$ given $\chi_n$ and $h^{-d/2}(\Lambda_n - k)$. As (6.5) implies that $1 - \alpha = P(h^{-d/2}(\Lambda_n^* - k) < \hat{q}_{n\alpha}|\chi_n) = P(\hat{Z} < \hat{q}_{n\alpha}) + \mathrm{o}(1)$, it follows that $\hat{q}_{n\alpha} = \hat{z}_\alpha + \mathrm{o}(1)$ conditionally on $\chi_n$. A similar argument by using Theorem 3.1 leads to $q_{n\alpha} = z_\alpha + \mathrm{o}(1)$. As $\hat{\Sigma}(x) \xrightarrow{p} \Sigma(x)$ uniformly in $x \in S$, then $\sigma^2(K, \hat{\Sigma}) \xrightarrow{p} \sigma^2(K, \Sigma)$, and hence $\hat{z}_\alpha = z_\alpha + \mathrm{o}(1)$. Therefore, $\hat{q}_{n\alpha} = q_{n\alpha} + \mathrm{o}_p(1)$ and this completes the proof. $\square$

**Proof of Theorem 3.3.** It can be shown that Lemmas 6.2–6.6 continue to hold true when we work under the local alternative $H_{1n}$. In particular, (6.3) is still valid. By using a derivation that resembles very much that for obtaining (6.4), we have $\Lambda_n(\vec{h}) = \{\Lambda_{n11}(\vec{h}) + \Lambda_{n112}^a(\vec{h}) + \Lambda_{n113}^a(\vec{h})\}\{1 + \mathrm{O}_p(h^2)\} + \mathrm{o}_p(h^{d/2})$, where $\Lambda_{n11}(\vec{h})$ is defined in (6.4),

$$\Lambda_{n112}^a(\vec{h}) = n^{-1} h^d c_n \sum_{l,t}^{k} \int \sum_{i,j}^{n} K_{h_l}(x - X_i) K_{h_t}(x - X_j) \gamma_{lt}(x) \pi(x) f^{-1}(x)$$
$$\times \{\epsilon_{jt} \Gamma_{nl}(X_i) + \epsilon_{il} \Gamma_{nt}(X_j)\} \mathrm{d}x$$

and

$$\Lambda_{n113}^a(\vec{h}) = n^{-1} h^d c_n^2 \sum_{l,t}^{k} \int \sum_{i,j}^{n} K_{h_l}(x - X_i) K_{h_t}(x - X_j)$$
$$\times \gamma_{lt}(x) \Gamma_{nl}(X_i) \Gamma_{nt}(X_j) \pi(x) f^{-1}(x) \mathrm{d}x.$$

It can be shown that $E\{\Lambda_{n112}^a(\vec{h})\} = 0$ and that

$$\begin{aligned}
E\{\Lambda_{n113}^a(\vec{h})\} &= (n-1) h^d c_n^2 \int \Gamma_n^{\mathrm{T}}(x) V^{-1}(x) \Gamma_n(x) f^2(x) \pi(x) \mathrm{d}x \\
&\quad + c_n^2 \beta_l^{-d} \int \sum_{l,t}^{k} R(\beta_l/\beta_t) \Gamma_{nl}(x) \gamma_{lt}(x) \Gamma_{nt}(x) \pi(x) \mathrm{d}x \{1 + \mathrm{O}(h^2)\} \\
&= n h^d c_n^2 \int \Gamma_n^{\mathrm{T}}(x) V^{-1}(x) \Gamma_n(x) f^2(x) \pi(x) \mathrm{d}x + \mathrm{O}(c_n^2 + n h^{d+2} c_n^2) \\
&= h^{d/2} \beta(f, K, \Sigma, \Gamma) + \mathrm{O}(c_n^2 + n h^{d+2} c_n^2) + \mathrm{o}(h^{d/2}).
\end{aligned} \tag{6.6}$$

It is fairly easy to see that $\Lambda^a_{n112}(\vec{h}) = o_p(h^{d/2})$ and $\Lambda^a_{n113}(\vec{h}) = h^{d/2}\beta(f,K,\Sigma,\Gamma) + o_p(h^{d/2})$. From Lemma 6.6, $h^{-d/2}[\Lambda_{n11}(\vec{h}) - k] \xrightarrow{d} N(0, \sigma^2(K,\Sigma))$. The theorem now follows after combining these results. $\square$

**Proof of Theorem 3.4.** The proof follows directly from Theorem 3.3 and from the fact that $\hat{q}_{n\alpha} = q_{n\alpha} + o_p(1)$, which is established in the proof of Theorem 3.2. $\square$

## Acknowledgements

Both authors are grateful for financial support from the IAP research network nr. P5/24 and P6/03 of the Belgian Government (Belgian Science Policy) and National Science Foundation Grants SES-0518904 and DMS-0604563.